\documentclass[12pt,a4paper]{article}
\usepackage[latin2]{inputenc}
\usepackage[T1]{fontenc}
\usepackage[leqno]{amsmath}
\usepackage{amsfonts,amssymb,latexsym,exscale}
\usepackage{graphicx}
\usepackage{pgf}
\usepackage{xspace}
\usepackage{theorem}
\usepackage{hyperref}
\usepackage{floatflt}
\usepackage[absolute]{textpos}
\usepackage{mathptmx}
\usepackage{textcomp}

\theoremstyle{plain}

\newtheorem{corollary}{Corollary}
\newtheorem{lemma}{Lemma}
\newtheorem{proposition}{Proposition}

\begin{document}

\title{\LARGE{\bf{ A Gauss-Bonnet formula  for moduli spaces of Riemann surfaces}}}
\author{\Large{Enrico Leuzinger}}
\date{}


\maketitle

\begin{abstract} We prove a Gauss-Bonnet  theorem for (finite coverings of) moduli spaces of Riemann surfaces
endowed with the  McMullen metric. The proof  uses 
properties of  an exhaustion of moduli spaces by compact submanifolds with corners
and the  Gauss-Bonnet formula of Allendoerfer and Weil for Riemannian polyhedra. 
\end{abstract}


\maketitle

\section{Introduction}

Mapping class groups (acting on Teichm\"uller spaces)
and  arithmetic subgroups of semisimple Lie groups
(acting on symmetric spaces) are prominent objects investigated in geometric group theory. In the last four decades several authors
- beginning with W. Harvey \cite{Har1}, J. Harer \cite{Ha} and N. Ivanov \cite{Iv1} - have discovered   remarkable similarities (but also differences) between these 
two classes of groups (see e.g. \cite{Iv3}, Section 9 for a discussion and also the more recent  \cite{Ar}, \cite{FLM} and \cite{Ji}).

The chief goal of the present paper is to add a further item to the list of similarities: the Gauss-Bonnet formula.
We consider a compact, orientable surface of genus $g$ with $p$ punctures  $S=S_{g,p}$ such that $3g-3+p>0$. The mapping class
group $\textup{Mod}(S)$ of $S$ acts properly discontinuously on Teichm\"uller space $\mathcal T(S)$.  Let $\Gamma$ be a torsion-free, finite index subgroup
of $\textup{Mod}(S)$. Then the Euler-Poincar\'e characteristic of $\textup{Mod}(S)$ is defined as 
$$
\chi(\textup{Mod}(S)):=[\textup{Mod}(S):\Gamma]^{-1}\ \chi(\Gamma\backslash\mathcal T(S)).
$$
Equivalently $\chi(\textup{Mod}(S))$ can be interpreted as the orbifold Euler characteristic of the moduli space of Riemann
surfaces $\textup{Mod}(S)\backslash\mathcal T(S)$.
J. Harer and D. Zagier \cite{HZ} showed that for $g>1$
$$
\chi(\textup{Mod}(S_{g,1}))=\zeta(1-2g)\ \ \ \textup{and also}\ \ \ \chi(\textup{Mod}(S_{g,0}))=\frac{1}{2-2g}\zeta(1-2g),
$$
where $\zeta(s)$ is the Riemann zeta function. Note that $\zeta(1-2g)=-B_{2g}/2g\in \mathbb Q$ where $B_{2g}$ is the $2g$-th Bernoulli number.

This result should be compared with the Euler-Poincar\'e characteristic of arithmetic groups. For the integer symplectic group, for instance,
a general result of G. Harder implies
$$\chi(Sp(2n, \mathbb Z))=\prod_{k=1}^n\zeta(1-2k)=\prod_{k=1}^n(-B_{2k}/2k),
$$
see \cite{H} or \cite{Se}, Section 3.
 Harder proved a Gauss-Bonnet theorem by which the computation of the Euler characteristic of an arithmetic subgroup $\Lambda$ of a semisimple Lie group $G$ becomes equivalent to the computation of the volume of the
 locally symmetric space  $\Lambda\backslash G/K$. The proof strongly relies on the reduction theory of arithmetic groups due to A. Borel and Harish-Chandra. In \cite{L4} I developped an analogous theory for mapping class groups based on work of Ivanov.
  This raises the question whether a Gauss-Bonnet theorem also holds in the setting of mapping class groups. In contrast to the symmetric case,
 where there is an essentially unique left-invariant metric,
 there is no canonical metric on Teichm\"uller resp. moduli space. However, as we will see,  the K\"ahler hyperbolic complete metric constructed by C. McMullen in \cite{Mc1} has 
 sufficiently nice properties such that a Gauss-Bonnet formula holds.

The generalized 
 Gauss-Bonnet theorem of
C.B.\ Allendoerfer, A.\ Weil and S.S.\ Chern  asserts that the Euler characteristic of
 a {\it closed} Riemannian manifold $(M, g)$ of even dimension $n$ is
given by
$$ \chi(M) = \int_{M}\omega_g = \int_{M}\Psi_g  \textup{dvol}_g$$
where the Gauss-Bonnet-Chern form $\omega_g = \Psi_g \, \textup{dvol}_g$  is (locally) computable
from the metric $g$. More precisely, let $\{\frac{\partial}{\partial y_k}\mid_p; 1\leq k\leq n\}$ be an orthonormal frame field in a neighbourhood of  $p\in M$
and let
$$
R_{ijkl}(p):=g( R( \frac{\partial}{\partial y_k}\mid_p,\frac{\partial}{\partial y_l}\mid_p)\frac{\partial}{\partial y_j}\mid_p,\frac{\partial}{\partial y_i}\mid_p),
$$
be the components of the curvature tensor, then
$$
\Psi_g(p)=\frac{1}{(2\pi)^{n/2}2^n(n/2)!}\sum_{\mu,\nu}\textup{sign}(\mu)\textup{sign}(\nu)R_{\mu_1\mu_2\nu_1\nu_2}(p)
R_{\mu_3\mu_4\nu_3\nu_4}(p)\ldots R_{\mu_{n-1}\mu_n\nu_{n-1}\nu_n}(p),
$$
where $\mu$ and $\nu$ independently run through all permutations of $\{1,2,\ldots, n\}$,
see \cite{AW}, \cite{C} or \cite{Spi}, Chapter 13 for more details.

Note that moduli spaces are even dimensional but not compact. However,  with respect to the McMullen metric they have bounded sectional curvature
 and finite volume (see \cite{Mc1}).
 The  Gauss-Bonnet theorem for
 {\it open} complete
  Riemannian manifolds with bounded sectional curvature
 and finite volume has been investigated by J.\ Cheeger and M.\ Gromov. 
They in particular showed  that such manifolds $M$ admit 
 exhaustions   by compact manifolds with smooth boundary, $M_i$, such that 
 $\textup{Vol} (\partial M_i)\rightarrow 0\ \  (i\rightarrow
\infty)$ and for which the second  fundamental forms $\textup{II}(\partial M_i)$
are uniformly bounded (see \cite{CG1}, \cite{CG2}, \cite{CG3} and also \cite{G} 4.5.C). By  a formula
 of Chern one has $ \chi(M_i) = \int_{M_i}\omega_g + \int_{\partial
M_i}\eta_i$
 where $\eta_i$ is a certain form on the boundary
$\partial M_i$ (see \cite{C}).  The above two properties  imply   that 
$\lim_{i\rightarrow\infty} \int_{\partial M_i}\eta_i = 0 $ 
and hence  that
$ \chi(M_i) =  \int_{M}\omega_g $  for sufficiently large $i$. As a consequence   the Gauss-Bonnet theorem holds {\it provided}
  $\chi(M_i) = \chi(M)$ for all sufficiently large $i$. This last property is not necessarily true for the exhaustions
  guaranteed by \cite{CG1}. It follows however, if in addition $\pi_1(M)$ is residually finite,  the injectivity radius of the universal
   covering space has a positive lower bound
  and  $M$ is of finite topological type (see \cite{CG1}, Thm. 3.1).

Locally symmetric spaces of higher rank  yield important  examples for the above class 
considered by Cheeger-Gromov.
For such (arithmetic) quotients the  Gauss-Bonnet formula was 
first  proved  by  Harder (see \cite{H}). In \cite{L1} I gave a new, simpler proof based on an 
  exhaustion $V = \bigcup_{s\geq 0}V(s)$ of locally symmetric spaces  by
{\it polyhedra}, i.e. compact submanifolds
with corners, and such that
 each polyhedron $V(s)$ is a strong deformation retract of $V$ (see \cite{L2}).
 The
 essential new feature of this exhaustion, which simplified the proof, is 
that the boundaries of 
 $\partial V(s)$  consist of subpolyhedra of $V(s)$ which are projections of pieces of horospheres in the universal covering space of $V$.
As a consequence 
 their second fundamental forms are uniformly bounded and a formula of Allendoerfer-Weil from the 1940s applies.

Coming back to moduli spaces of Riemann surfaces, we
note that the (orbifold-)
fundamental group of moduli space is  the mapping class group, which is known to be residually finite (see \cite{Iv2}). Moreover, with respect to the McMullen metric,  the injectivity radius 
of the universal covering space (i.e., Teichm\"uller space) has a positive lower bound (see \cite{Mc1}) and finally moduli space is of finite topological type (see e.g. Lemma 1 below). Thus by the result of Cheeger-Gromov refered to above, a Gauss-Bonnet formula holds. The goal of the present note is to give a new, independent proof of this fact.
In contrast to the general situations considered by Cheeger and Gromov there is additional structural information available (as in the in the case of locally symmetric spaces).
This allows us to avoid the rather involved technical arguments and constructions in \cite{CG2}.

\vspace{2ex}

\noindent{\bf Theorem} \ {\bf (Gauss-Bonnet for moduli spaces)} {\it Let $S=S_{g,p}$ be a compact, orientable surface of genus $g$ with $p$ punctures such that $3g-3+p>0$.
 Let $\Gamma\subset \textup{Mod}(S)$ be   a torsion free, finite index subgroup of the mapping class group of $S$ consisting of pure elements and let  $\mathcal M(S)= 
\Gamma\backslash\mathcal T(S)$ be the corresponding finite covering space of the moduli space of 
Riemann surfaces $\textup{Mod}(S)\backslash\mathcal T(S)$.  Let $\omega_{g_M}=\Psi_{g_M}  \textup{dvol}_{g_M} $ be the Gauss-Bonnet-Chern form associated to the   McMullen metric $g_M$ on $\mathcal M(S)$
and let $\chi(\mathcal M(S))$ be the Euler-Poincar\'e characteristic of $\mathcal M(S)$. Then the following
 Gauss-Bonnet formula holds: 
$$\chi(\mathcal M(S)) = \int_{\mathcal M(S)}\omega_{g_M}= \int_{\mathcal M(S)}\Psi_{g_M}   \textup{dvol}_{g_M}.$$
}

The proof is modelled on the one for locally symmetric spaces in \cite{L1}. We  use an exhaustion of moduli spaces by compact submanifolds with corners due to Ivanov \cite{Iv1}
whose properties were further  studied in \cite{L4}. As already mentioned above,  moduli space has 
  bounded sectional curvature and finite volume with respect to the McMullen metric.  We will show in addition that the boundaries of the exhausting polyhedra  have bounded second fundamental forms, so that we again can use \cite{AW}.

\vspace{3ex}

{\bf Remarks}. (1)\  There are several metrics on Teichm\"uller resp. moduli space which are quasi-isometric to the Teichm\"uller metric
(see e.g. \cite{LSY}). It is conceivable that
a  Gauss-Bonnet formula also  holds for certain of these other metrics (some of which seem to be more canonical from the point of view of possible applications). 

(2) In principle, uniform bounds on the density $\Psi_{g_M}$ in the above Gauss-Bonnet formula can be used to estimate volumes of moduli spaces with respect to the McMullen (or Teichm\"uller)
metric. Unfortunately, mere uniform curvature bounds  $|R_{ijkl}(p)|\leq \textup{const.}$ seem to be too coarse to provide interesting asymptotic 
volume information. Note that for the surfaces $S_{g,0}$ the (Teichm\"uller) volume grows  more than exponentially in $g$. In fact, Proposition 2 below and the remark following it together with Proposition 3 imply that  it is at least $g!$. This is similar to the volume growth of locally symmetric space with respect to dimension.

\vspace{3ex}

{\bf Notation}.\ Explicit constants are irrelevant for our purpose.  If $f$ and $g$ are positive real valued functions 
on a set $S$
we thus simply write $f\prec g$ if  there is a constant $c > 0$ such that $f(s)\leq c g(s)$ for all $s\in S$. Similarly, $f\asymp g$ if
$f(s)\leq c_1 g(s)$ and $g(s)\leq c_2 f(s)$ for positive constants
$c_1,c_2$.
\vspace{3ex}

{\bf Acknowledgement}. \ I would like thank C.T. McMullen  for helpful correspondence.

\section{The formula of Allendoerfer and Weil}

A  $C^{\infty}$-manifold with corners is a topological Hausdorff space
 locally modeled 
upon a product of lines and half-lines and such that  coordinate changes
 are  of class $C^{\infty}$. For precise definitions and basic
 information about
this concept  we refer to  \cite{DH}.
A {\it Riemannian polyhedron} is a compact manifold with corners equipped with a Riemannian metric.

Let $P^n$  be an $n$-dimensional Riemannian polyhedron with  boundary consisting of a finite family of lower dimensional subpolyhedra
$P_E^{n-k}$  ($1\leq k\leq n$) and with Riemannian metric induced from
$ P^n$. The {\it outer angle} $O(p)$
 at a point $p$ of $P_E^{n-k}$
 is defined as the set  of all unit tangent vectors $v\in T_p  P^n$ such that
$\langle v, w\rangle_p\leq 0$ for all $w$ in the tangent cone of $ P^n$ at $p$.
Note that $O(p)$ is a spherical cell bounded by ``great spheres''
in the $(k-1)$-dimensional unit sphere of the normal space of $ P_E^{n-k}\subset  P^n$  at $p$.
In \cite{AW}
Allendoerfer and Weil define a certain real valued 
function $\Psi_{E,k}$ on the outer angles of $ P_E^{n-k}$.
The explicit form of this function will not be needed in this paper. We shall only use the fact 
that   $\Psi_{E,k}$ is
 locally computable from the components of
the metric and the curvature tensor of $P^n$ and from the components of the second
fundamental forms $\textup{II}_Z(p), Z\in O(p)$, of $ P_E^{n-k}$ in $ P^n$.
 Let $\Psi\, \textup{dvol}$ denote the Gauss-Bonnet-Chern form on 
$ P^n$ and  $\textup{dvol}_E$ (resp. $d\omega_{k-1}$) the volume element
of $ P_E^{k}$ (resp. of the standard unit sphere $S^{k-1}$).
The {\it inner Euler characteristic} $\chi'$ of $ P^n$  is by definition the
Euler characteristic of the open complex consisting of all inner cells in an
arbitrary simplicial subdivision of $ P^n$.

The generalized Gauss-Bonnet
formula of Allendoerfer-Weil  for Riemannian polyhedra proved in \cite{AW} has  then the following form:

\begin{proposition}\label{aw}
Let $ P^n$ be a Riemannian polyhedron of even dimension $n$  with boundary consisting of a finite family of subpolyhedra 
$ P_E^{n-k}$. Then the inner Euler characteristic of $ P^n$
is given by $$ \chi'( P^n) = \int_{ P^n}\Psi \, \textup{dvol} + \sum_{E}\sum_{k=1}^{n}
\int_{ P_E^{n-k}}\ (\int_{O(p)}\Psi_{E,k} \ d\omega_{k-1})  \ \textup{dvol}_E(p).$$
\end{proposition}

\vspace{3ex}

{\bf Remark}.\ This formula has  also been used recently by McMullen in \cite{Mc2} to compute the complex hyperbolic volume
of moduli spaces of configurations of  ordered points on the Riemann sphere.

\section{A polyhedral exhaustion of moduli spaces }

\subsection{Moduli spaces and the complex of curves}

Let $S=S_{g,p}$ be a compact, orientable surface  of genus $g$ with $p$ punctures
such that $d(S):=3g-3+p>0$. This last assumption implies that $S$ carries  finite volume Riemannian metrics of constant  curvature $-1$ and $p$  cusps. 
A {\it marked hyperbolic surface} is a pair $(X, [f])$ where $X$ is a smooth surface equipped with a complete Riemannian metric of constant curvature $-1$ and where $[f]$ denotes the isotopy class of  a diffeomorphism   $f:X\longrightarrow S$
mapping cusps to punctures.
Two marked surfaces 
$(X_1,[f_1])$ and
$(X_2,[f_2])$ are equivalent
if there is an isometry  $h:X_1\longrightarrow X_2$ such that
$[f_2\circ h]=[f_1]$.  The collection of these equivalence classes is (one possible definition of) \emph{the Teichm\"uller space}
$\mathcal{T}(S)$ of $S$.
 The
 {\it mapping class group} $\textup{Mod}(S)$ is the group of all orientation preserving 
 diffeomorphisms of $S$, which fix the punctures,
modulo isotopies which also fix the punctures.  The group $\textup{Mod}(S)$ acts on $\mathcal T(S)$ 
according to the formula $h\cdot (X,[f]):=(X,[h\circ f])$, for $ h\in \textup{Mod}(S), (X,[f])\in \mathcal T(S)$.
 The corresponding orbit space $\textup{Mod}(S)\backslash \mathcal T(S)$ is the {\it moduli space} of isometry classes of hyperbolic surfaces 
  (obtained by forgetting the marking).

The 
\emph{complex of curves} $\mathcal C(S)$ of $S$ is an infinite (even locally infinite) simplicial complex
of dimension $d(S)-1$. Note that $d(S)=\frac{1}{2}\dim \mathcal T(S)$.
The \emph{vertices} of $\mathcal C(S)$  are the isotopy classes of simple closed curves (called \emph{circles}) on $S$, which are non-trivial (i.e., not contractible in $S$ to a point or to a component of $\partial S$). We denote the isotopy class
of a circle $C$ by $\langle C\rangle$. A set of $k+1$ vertices $\{\alpha_0,\ldots \alpha_k\}$ spans a
$k$-\emph{simplex} of $\mathcal C(S)$ if and only if $\alpha_0=\langle C_0\rangle,\ldots,\alpha_k=\langle C_k\rangle$
for some set of pairwise non-intersecting circles $C_0,\ldots, C_k$. 
For a simplex $\sigma\in \mathcal C(S)$ we denote by $|\sigma|$ the number of its vertices. A crucial fact  is that $\mathcal C(S)$ is a thick chamber complex. In particular every simplex is the face of a maximal simplex.  Moreover the mapping class group $\textup{Mod}(S)$ acts simplicially on $\mathcal C(S)$ and the quotient  $ \textup{Mod}(S)\backslash\mathcal C(S)$ is a finite (orbi-) complex  (see \cite{Har1}, Proposition 1).

We next consider finite index, torsion-free subgroups $\Gamma$ of $\textup{Mod}(S)$ which, in addition, consist of pure mapping classes. Recall that a mapping class 
$h\in \textup{Mod}(S)$ is called {\it pure} if it can be represented by a diffeomorphism $f:S\longrightarrow S$ fixing (pointwise) some union $\Lambda$ of disjoint and pairwise non-isotopic nontrivial circles on $S$ and such that $f$ does not permute the components of $S\setminus \Lambda$ and induces on each component of the cut surface $S_{\Lambda}$ a diffeomorphism isotopic to a pseudo-Anosov
or to the identity diffeomorphism (see \cite{Iv3}, \S 7.1). It is well-known that such subgroups 
$\Gamma$ exist. For example one can take $\Gamma=\Gamma_S(m)$, the kernel of the natural homomorphism $\textup{Mod}(S)\longrightarrow \textup{Aut}(H_1(S,\mathbb Z/m\mathbb Z)), m\geq 3$,
defined by the action of diffeomorphisms on homology (see e.g. \cite{Iv3}, \S 7.1).

\subsection{A partition and exhaustion of moduli space}

Let $\Gamma \subset \textup{Mod}(S)$ be as in the previous section. Using work of N. Ivanov I constructed in 
 \cite{L4} 
a  $\Gamma$-invariant partition of Teichm\"uller space into disjoint subsets, which in turn yields a
partition of   $\mathcal M(S)= \Gamma\backslash \mathcal T(S)$.
In order to formulate that result more precisely, we need  {\it length functions}.
Let $\alpha$ be a vertex of $\mathcal C(S)$, i.e., $\alpha=\langle C\rangle$
for a (non-trivial) circle $C$. 
Since $d(S)>0$,  any point $X\in \mathcal T(S)$ represents a finite volume Riemann surface of  curvature $-1$ with $p$ cusps. On the surface $X$ the isotopy class of $\alpha$ 
contains a unique geodesic loop; let $l_{\alpha}(X)$ denote its length. This defines 
a (smooth) function $l_{\alpha}: \mathcal T(S)\longrightarrow \mathbb R_{>0}$ for every vertex $\alpha\in \mathcal C(S)$. 
For $\varepsilon>0$ we then define the {\it $\varepsilon$-thick part of Teichm\"uller space} 
$$\textup{Thick}_{\varepsilon}\mathcal T(S):=
\{X\in\mathcal T(S)\mid l_{\alpha}(X)\geq \varepsilon, \forall \alpha\in\mathcal C(S)\}.
$$
 This set is 
$\Gamma$-invariant and its quotient, $\textup{Thick}_{\varepsilon}\mathcal M(S)$, is the {\it $\varepsilon$-thick part  of $\mathcal M(S)= \Gamma\backslash \mathcal T(S)$}.

For every vertex $\alpha\in \mathcal C(S)$ we further set $\mathcal H_{\varepsilon}(\alpha):=\{X\in \textup{Thick}_{\varepsilon}\mathcal T(S)\mid l_{\alpha}(X)=\varepsilon\}$.
Then (for fixed $\varepsilon$ sufficiently small) we have $\partial \textup{Thick}_{\varepsilon}\mathcal T(S)=\bigcup_{\alpha} \mathcal H_{\varepsilon}(\alpha)$. Given $\sigma \in \mathcal C(S)$ we denote by $S_{\sigma}$ the corresponding  cut surface, i.e., the 
result of cutting $S$ along (non-intersecting) circles from the isotopy classes  $\alpha\in\sigma$.
In \cite{L4} I proved 

\begin{proposition}\label{partitionmod}Let $S=S_{g,p}$ be a compact, orientable surface of genus $g$ with $p$ punctures such that $3g-3+p>0$. Let $\Gamma$ be   a torsion free, finite index subgroup of the mapping class group of $S$ consisting of pure elements and let  $\mathcal M(S)= 
\Gamma\backslash\mathcal T(S)$ be the corresponding finite covering space of the moduli space of 
Riemann surfaces. Finally  let $\mathcal E$ be the set of simplices of the  finite simplicial complex $ \Gamma\backslash\mathcal C(S)$. 

Then following holds:

$\textup{(1)}$\  There exists $\varepsilon =\varepsilon (S)>0$ such that there is a partition of  $\mathcal M(S)$, i.e.,  a \textup{disjoint union}
$$
\mathcal M(S)= \textup{Thick}_{\varepsilon}\mathcal M(S)\sqcup\bigsqcup_{\sigma\in \mathcal E}\textup{Thin}_{\varepsilon}(\mathcal M(S),\sigma),
$$ 
where $\textup{Thick}_{\varepsilon}\mathcal M(S)$ is a compact submanifold with corners  and 
each  $\textup{Thin}_{\varepsilon}(\mathcal M(S),\sigma)$ is diffeomorphic to
 $\mathcal B_{\varepsilon}(\sigma)\times \mathbb R_{>0}^ {|\sigma|}$ with 
$\mathcal B_{\varepsilon}(\sigma)$  a (trivial) torus bundle over the $\varepsilon$-thick part of  moduli space of the cut surface
$S_{\sigma}$ such
that the length  of each boundary circle is $\varepsilon$:
$$
0\rightarrow T_{\varepsilon}^{|\sigma|}\rightarrow \mathcal B_{\varepsilon}(\sigma)
\rightarrow \textup{Thick}_{\varepsilon}\mathcal M(S_{\sigma})\rightarrow 0.
$$

$\textup{(2)}$\  Each boundary face $\mathcal B_{\varepsilon}(\sigma)$, for $\sigma\in \mathcal E$, of $\textup{Thick}_{\varepsilon}\mathcal M(S)$ is the image of $\bigcap_{\alpha\in\sigma}\mathcal H_{\varepsilon}(\alpha)$
under the natural projection $\pi:\mathcal T(S)\longrightarrow \mathcal M(S)$.
\end{proposition}

\vspace{3ex} 

{\bf Remark}.\ In a similar exhaustion for locally symmetric spaces, the thin parts are indexed by the simplices
 of the rational Tits building modulo the arithmetic group (see \cite{L3}, Theorem C). The number of cells in this (finite) quotient is rather small; for instance, for the locally symmetric space of 
 principally polarized abelian varieties, $Sp(2n,\mathbb Z)\backslash Sp(2n,\mathbb R)/U(n)$, it is always one. In contrast, the corresponding quotient for moduli spaces is large: The number of maximal simplices in $ \textup{Mod}(S)\backslash\mathcal C(S)$ equals the number of different pant decompositions of $S$, which is roughly $g!$ for large $g$ (see \cite{Bu}, Theorem 3.5.3).

\vspace{3ex}

Proposition 2 yields  the

\begin{corollary}  \label{exhaustion} There exists  $\varepsilon_0>0$ depending only on $S$, such that there is a $\Gamma$-invariant exhaustion of Teichm\"uller space $\mathcal T (S) = \bigcup_{\varepsilon\leq \varepsilon_0}\textup{Thick}_{\varepsilon}\mathcal T(S)$, which induces an exhaustion of 
$\mathcal M(S)$  by  polyhedra, i.e., by compact submanifolds  with corners: $\mathcal M (S) = \bigcup_{\varepsilon\leq \varepsilon_0}\textup{Thick}_{\varepsilon}\mathcal M(S)$.
\end{corollary}

In contrast to some of the exaustions constructed by Cheeger-Gromov (see Section 1), the topology of the exhausting polyhedra does not change.
Namely, we have

\begin{lemma}\label{retraction} For every sufficiently small $\varepsilon$ the
 moduli space $\mathcal M(S)$ is homeomorphic to the interior of the
polyhedron $\textup{Thick}_{\varepsilon}\mathcal M(S)$ in $\mathcal M(S)$, and  $\textup{Thick}_{\varepsilon}\mathcal M(S)$  is 
a strong deformation retract of $\mathcal M(S)$.
\end{lemma}

\emph{Proof}.\ 
Recall that the complex of curves $\mathcal C(S)$ is a thick chamber complex.
Given a simplex $\sigma\in \mathcal C(S)$ we can thus choose a simplex $\tau\in \mathcal C(S)$ of maximal dimension 
$d(S)-1=3g-4+p$ containing $\sigma$. Then, by \cite{Abi},  there are  \emph{adapted Fenchel-Nielsen coordinates} 
on Teichm\"uller space $\mathcal T(S)$, i.e., there is  a diffeomorphism
$$
\Phi_{\tau}: \mathcal T(S)\longrightarrow \mathcal T(S_{\sigma})\times \mathbb R^{|\sigma|}\times \mathbb R_{ > 0}^{|\sigma|}\ ;\ X\longmapsto (s(X), \theta(X),l(X)),
$$
where $s=(\theta_{\beta},l_{\beta})_{\beta\in\tau\setminus\sigma}$ parametrizes the Teichm\"uller space $\mathcal T(S_{\sigma})$ of the cut surface $S_{\sigma}$, 
$\theta=(\theta_{\alpha})_{\alpha\in\sigma}\in \mathbb R^{|\sigma|}$ are twist parameters  and
$l= (l_{\alpha})_{\alpha\in\sigma}\in \mathbb R_{>
 0}^{|\sigma|}$ are length functions 
(here and elsewhere $|\sigma|$ denotes the number of vertices of $\sigma$). 

We define the \emph{outer cone } at $X\in\partial \textup{Thick}_{\varepsilon}\mathcal T(S)$  as the preimage 
$$CO(X):=\Phi_{\tau}^{-1}\{(s(X),  \theta(X),(l_{\alpha})_{\alpha\in \sigma})\mid l_{\alpha}    <     \varepsilon\ \ \textup{for all}\ \      \alpha\in \sigma\}.
$$

By \cite{L4}, Lemma 5, outer cones are $\Gamma$-invariant. This
allows one to define {\it outer cones in moduli space}\,: for $X\in \textup{Thick}_{\varepsilon}\mathcal M(S)$
we set $CO(X):= \pi(CO(\hat{X}))$ where $\hat{X}$ is any lift of $X$. We then have (see \cite{L4}, 2.2):
$$
\mathcal M(S)=\textup{Thick}_{\varepsilon}\mathcal M(S)\sqcup \bigsqcup_{X\in\partial\textup{Thick}_{\varepsilon}\mathcal M(S)}CO(X)
= \textup{Thick}_{\varepsilon}\mathcal M(S)\sqcup \bigsqcup_{\sigma\in \mathcal E}\  
\bigsqcup_{X\in\mathcal B_{\varepsilon}(\sigma)}CO(X).
$$
Note that each outer cone  is diffeomorphic to the open {\it hyperoctant} $\mathbb R_{>0}^{|\sigma|}$. 
The claimed retraction is then  given by retracting each outer cone  to the apex of its closure.
\hfill$\Box$

\section{Estimates for the boundary subpolyhedra}

We wish to apply Proposition \ref{aw} to the polyhedra $\textup{Thick}_{\varepsilon}\mathcal M(S)$ in the above exhaustion
and then take the limit for $\varepsilon\to 0$.
To that end we need uniform estimates for the second fundamental forms and the volumes  of the
(lower dimensional) boundary polyhedra with respect to the McMullen metric. 

By Section 3, Proposition \ref{partitionmod}, we know that for every (sufficiently small) $\varepsilon$ the boundary polyhedra are
of the form $\mathcal B_{\varepsilon}(\sigma)$ for $\sigma\in \mathcal E$, the set of simplices in the finite simplicial complex
$ \Gamma\backslash\mathcal C(S)$. Moreover, if $n:=\dim \mathcal M(S)$, we have $\dim \mathcal B_{\varepsilon}(\sigma)=n-
|\sigma|$.

\subsection{An approximation of the metric on thin parts}

We will prove the Gauss-Bonnet formula with  respect to the metric on moduli space constructed by McMullen 
in \cite{Mc1}.
He  showed that with respect to this metric $\mathcal M(S)$  
is K\"ahler hyperbolic in the sense of Gromov and thus in particular carries a complete
 finite volume 
 Riemannian metric of
  bounded sectional curvature. The McMullen metric is a modification of  the (incomplete) 
Weil-Petersson metric 
and quasi-isometric to the Teichm\"uller metric (see \cite{Mc1}).

We next describe explicit representatives of the bi-Lipschitz class of the Teichm\"uller resp. McMullen metric
on $(\varepsilon,\sigma)$-thin parts of  $\mathcal T(S)$ (and $\mathcal M(S)$). To that end we use
 an  expansion of the Weil-Petersson (WP) metric due to  S. Wolpert 
(see \cite{Wo1}, \cite{Wo2}). 

For each  length function $l_{\alpha}$ we set  $u_{\alpha}:=-\log l_{\alpha}^{1/2}$. 
Considering this logarithmic root
 length instead of $l_{\alpha}$ itself is suggested by  work of  Wolpert (see e.g. 
\cite{Wo1}, \cite{Wo2}). 
Following Wolpert we also set $\lambda_{\alpha}:= \textup{grad}\  l_{\alpha}^{1/2}$
 (resp. 
$\nu_{\alpha}:= -\textup{grad}\  u_{\alpha}$) and define the  \emph{Fenchel-Nielsen-gauge} 
as the differential 1-form 
$\rho_{\alpha}:=
2\pi(l_{\alpha}^{3/2}\langle\lambda_{\alpha},\lambda_{\alpha}\rangle)^{-1}\langle\ ,
J\lambda_{\alpha}\rangle$
$=2\pi(\langle\nu_{\alpha},\nu_{\alpha}\rangle)^{-1}\langle\ ,J\nu_{\alpha}\rangle$ 
(see \cite{Wo2}, 4.15). This gauge is
 normalized such that $l_{\alpha}(T_{\alpha})=1$ for $T_{\alpha}:=
(2\pi)^{-1}l_{\alpha}^{3/2}J\lambda_{\alpha}$  the WP-unit infinitesimal FN angle variation. We also set
$\tilde{\rho}_{\alpha}:=l_{\alpha}^{-1/2}\rho_{\alpha}$ (this gauge is normalized  with respect to the Teichm\"uller-unit infinitesimal FN angle variation).
In  \cite{L4} I showed

\begin{proposition} \label{approxmetric}There is $\varepsilon_*>0$ depending  only on $S$,
 such that for $\sigma\in\mathcal C(S)$, 
 $\varepsilon \leq\varepsilon_*$  and the   $(\varepsilon,\sigma)$-thin part 
$$
\textup{Thin}_{\varepsilon}(\mathcal T(S),\sigma):=
\bigsqcup_{X\in\bigcap_{\alpha\in\sigma}\mathcal H_{\varepsilon}(\alpha)} CO(X) 
\cong\textup{Thick}_{\varepsilon}\mathcal T(S_{\sigma})\times 
\mathbb R^{|\sigma|}\times \mathbb R_{> 0}^{|\sigma|}
$$ 
the  following Finsler metric expansion of the Teichm\"uller (or McMullen) metric with respect to
 adapted FN-coordinates
 holds 
$$
\|.\|^2_{\mathcal T(S)}\asymp 
\|.\|^2_{\mathcal T(S_{\sigma})}
+ \sum_{\alpha\in \sigma}e^{-6u_{\alpha}}\tilde{\rho}_{\alpha}^2+du^2_{\alpha}.
$$
The bi-Lipschitz constants involved in this estimate only depend on $S$ and $\varepsilon_*$.
\end{proposition}

This metric approximation decends to the $(\varepsilon,\sigma)$-thin parts of moduli space (see \cite{L4}, Cor. 3):

$$
\textup{Thin}_{\varepsilon}(\mathcal M(S),\sigma):=
\bigsqcup_{X\in \mathcal B_{\varepsilon}(\sigma)} CO(X) 
\cong\textup{Thick}_{\varepsilon}\mathcal M(S_{\sigma})\times 
T_{\varepsilon}^{|\sigma|}\times \mathbb R_{> 0}^{|\sigma|}.
$$ 

\subsection{Second fundamental forms of boundary polyhedra}

\begin{lemma}\label{2ndff} Let $\varepsilon_0$ be so small that both Proposition \ref{partitionmod} and 
Proposition \ref{approxmetric} hold. Take $\varepsilon\leq \varepsilon_0$ and $\sigma\in \mathcal E$.
 Then the second fundamental forms of
 the boundary polyhedron $\mathcal B_{\varepsilon}(\sigma)$ with respect to outer angles in  $\textup{Thick}_{\varepsilon}\mathcal M(S)$ are uniformly bounded
by a constant independent of $\sigma$ and $\varepsilon$.
\end{lemma}

\emph{Proof}.\ By Proposition \ref{partitionmod}(2) the boundary polyhedra
are projections of intersections of level sets of length functions
$$\mathcal B_{\varepsilon}(\sigma)=\pi(\bigcap_{\alpha\in\sigma}\mathcal H_{\varepsilon}(\alpha))\ \ \ \textup{where \ }\ \ \mathcal H_{\varepsilon}(\alpha):=\{X\in \textup{Thick}_{\varepsilon}\mathcal T(S)\mid l_{\alpha}(X)=\varepsilon\}.
$$
The outer angles (see Section 2) of $\mathcal B_{\varepsilon}(\sigma)$ are thus positive linear combinations of gradients of the lenght functions $ l_{\alpha}$ for $\alpha\in \sigma$. 
Set $N_{\alpha}:=\frac{\textup{grad}(l_{\alpha})}{\|\textup{grad}(l_{\alpha})\|}$. Hence, to prove the Lemma, we have  to show that (with respect to the McMullen metric)
$$\langle D_XY,N_{\alpha}\rangle=O(1)\ \ 
\textup{for all}\ \  \alpha\in \sigma \ \ \ \textup{and all  unit tangent vectors }\ \  X,Y\ \ \textup{of}\ \ \mathcal B_{\varepsilon}(\sigma).
$$

We first note that  the  derivatives up to order two of $l_{\alpha}, \alpha\in \sigma,$  are $O(l_{\alpha})$.
In fact, in the proof of Theorem 8.2. in \cite{Mc1} McMullen used the Bers embedding and euclidean polydiscs to show that the derivatives
of $\log(l_{\alpha})$ are $O(1)$ and hence the first dervatives of  $l_{\alpha}$ are $O(l_{\alpha})$ and similarly for second derivatives.

We next show that $ \|\textup{grad}(l_{\alpha})\|\asymp l_{\alpha}$. 
Denote by $\textup{grad}_*(l_{\alpha})$ the gradient of $l_{\alpha}$
with respect to the approximating product metric in Proposition  \ref{approxmetric}. Then for the adopted FN-coordinates we have $l_{\alpha}=e^{-2u_{\alpha}}$ hence, if we denote by $\langle \cdot\ ,\cdot \rangle_*$  the approximate metric, we have that $\textup{grad}_*(l_{\alpha})$ is parallel to $\partial_{u_{\alpha}}$
and $\| \partial_{u_{\alpha}} \|_*=1$.
Hence
$$
\|\textup{grad}_*(l_{\alpha})\|_*=\langle\textup{grad}_*(l_{\alpha}), \partial_{u_{\alpha}}\rangle_*=dl_{\alpha}(\partial_{u_{\alpha}})=-2l_{\alpha}.
$$
Since (on the
thin part $\textup{Thin}_{\varepsilon}(\mathcal M(S),\sigma)$ under consideration) the McMullen metric is comparable to the approximate metric  
we eventually find 
$$
\|\textup{grad}(l_{\alpha})\|\asymp \|\textup{grad}_*(l_{\alpha})\|_*\asymp
l_{\alpha}.
$$

By this estimate for the gradient,  
$\langle D_XY,N_{\alpha}\rangle$ is $O(1)$ if $\langle D_XY,\textup{grad}(l_{\alpha})\rangle$ is $O(l_{\alpha})$ for unit tangent vectors $X,Y$ of  $\mathcal B_{\varepsilon}(\sigma)$.
Now such an $Y$ is tangent to the level surface of $l_{\alpha}$ and thus $\langle Y,\textup{grad}(l_{\alpha})\rangle=0$. This yields (see \cite{ON},
Lemma 3.49)
$$
\langle D_XY,\textup{grad}(l_{\alpha})\rangle=-\langle Y,D_X\textup{grad}(l_{\alpha})\rangle=-\textup{Hess}_{l_{\alpha}}(X,Y)=-(XY)l_{\alpha}+(D_XY)l_{\alpha}=O(l_{\alpha})
$$
were the last estimate follows from the already mentioned estimates for the derivatives of order one and  two and the additional fact that the covariant derivative $D_XY$ is $O(1)$ for unit vectors $X,Y$  since the derivatives of the metric are also bounded (see again\cite{Mc1}, Theorem 8.2).

\hfill$\Box$

\subsection{Volumes of boundary polyhedra}

We next  estimate the volumes of the boundary polyhedra  $\mathcal B_\varepsilon(\sigma)$ which have dimension $n-|\sigma|$, where $n=\dim \mathcal M(S)$.

\begin{lemma}\label{finitevol}  The $(n-|\sigma|)$-dimensional volume of each 
boundary subpolyhedron  $\mathcal B_\varepsilon(\sigma)$ of $\textup{Thick}_{\varepsilon}\mathcal M(S)$  satisfies
 $$\textup{Vol}\ (\mathcal B_\varepsilon(\sigma))\prec \varepsilon,$$
where the  constants involved in this estimate are independent of $\sigma$ and $\varepsilon$.
\end{lemma}

\emph{Proof}.\ Since boundedness of volume is a bi-Lipschitz invariant, we can work with the approximation of the McMullen metric from Proposition \ref{approxmetric}:
 On the thin part $\textup{Thin}_{\varepsilon}(\mathcal M(S),\sigma)$ the latter is comparable to a product of the McMullen metric on the thick part of the lower-dimensional moduli space $\mathcal M(S_{\sigma})$ with a product of 2-dimensional hyperbolic metrics (the number of factors of the latter beeing equal to  $|\sigma|$). Now,
 by the general properties of the McMullen metric,  the volume of $\mathcal M(S_{\sigma})$ is bounded, hence in particular
$\textup{Vol}(\textup{Thick}_{\varepsilon}\mathcal M(S_{\sigma}))\leq c_1(\sigma)$.
Similarly, we get for the torus fibre (see Proposition \ref{partitionmod} (1)) from the expansion in Proposition \ref{approxmetric}
$$\textup{Vol}(T_{\varepsilon}^{\sigma}))\leq c_2(\sigma)e^{-6(-\log\sqrt{\varepsilon}) |\sigma|}= c_2(\sigma) \varepsilon^{3|\sigma|}.$$
Thus
$$
\textup{Vol}(\mathcal B_{\varepsilon}^{\sigma}))\leq c_1(\sigma)c_2(\sigma) \varepsilon^{3|\sigma|}.$$

Since the set $\mathcal E$ of simplices indexing the boundary polyhedra is finite and $\varepsilon\ll 1$, the claim of the Lemma follows.
\hfill$\Box$

\section{The proof of the Gauss-Bonnet formula}

In this section we give the proof of the Gauss-Bonnet formula for moduli spaces of Riemann surfaces equipped with the McMullen metric
as stated in the introduction.

 By Corollary \ref{exhaustion} there is an exhaustion 
$\mathcal M (S) = \bigcup_{\varepsilon\leq \varepsilon_0}\textup{Thick}_{\varepsilon}\mathcal M(S)$ of $\mathcal M(S)$  by 
Riemannian polyhedra. Each  polyhedron  $\mathcal M(\varepsilon):=\textup{Thick}_{\varepsilon}\mathcal M(S)$ in this exhaustion 
 is equipped with the Riemannian metric induced by the one of $\mathcal M(S)$. We  set $n:=\dim \mathcal M(S)$ and $k:=|\sigma|$ (i.e.
  $\dim \mathcal B_{\varepsilon}(\sigma)=n-k$). Since $n$ is even
   Proposition \ref{aw} applied to $\mathcal M(\varepsilon)$ yields
$$ |\chi'(\mathcal M(\varepsilon))  - \int_{\mathcal M(\varepsilon)}\Psi\, \textup{dvol} | \prec \sum_{\sigma\in \mathcal E}\sum_{k=1}^{n}
\int_{\mathcal B_{\varepsilon}(\sigma)}\int_{O(X)}
\|\Psi_{\sigma,k}\|\  d\omega_{k-1} \ \ \textup{dvol}_{\sigma}(X).$$
 As remarked in Section 1,
the function  $\Psi_{\sigma,k}$ is locally computable from the components of
the metric and the curvature tensor of $\mathcal M(\varepsilon)$ and from the components of the second
fundamental forms of the $\mathcal B_{\varepsilon}(\sigma)$ in $\mathcal M(\varepsilon)$.
The fact that the McMullen metric has bounded curvature together with Lemma \ref{2ndff} thus implies    that  $\|\Psi_{\sigma,k}\|\prec 1$ for all
$\sigma, k$.  Using Lemma \ref{finitevol} we conclude that
 $$|\chi'(\mathcal M(\varepsilon))  - \int_{\mathcal M(\varepsilon)}\Psi\, \textup{dvol} |  \prec
  \sum_{k,\sigma}\textup{Vol}(\mathcal B_{\varepsilon}(\sigma))\prec \varepsilon.$$
By Lemma \ref{retraction} we have
 $\chi'(\mathcal M(\varepsilon)) = \chi(\mathcal M(S))$ and since  $\chi(\mathcal M(S))$ is an integer
 we have $\chi(\mathcal M(S)) = \int_{\mathcal M(\varepsilon)}\Psi\, \textup{dvol}$ for all sufficiently small $\varepsilon$.
 Since the polyhedra $\mathcal M(\varepsilon)$ exhaust $\mathcal M(S)$,   the claimed formula follows.\hfill$\Box$

\vspace{5ex}

\noindent\textsc{Institute for Algebra und Geometry\\
 Karlsruhe Institute of Technology (KIT),  Germany}

\vspace{1ex}

\noindent{\tt enrico.leuzinger@kit.edu}
\end{document}